\documentclass[12pt,twoside]{article}
\usepackage{amsmath,amssymb}
\usepackage{mathrsfs}
\usepackage{hyperref}
\usepackage{amsmath,amsthm}
\usepackage[mathscr]{euscript}

\newcommand{\N}{\mathbb{N}}
\newtheorem{theorem}{Theorem}
\newtheorem{lemma}{Lemma}

\newtheorem{definition}{Definition}
\begin{document}
\thispagestyle{empty}
\setcounter{page}{1}

\noindent
{\footnotesize {\rm To appear in\\[-1.00mm]
{\em Dynamics of Continuous, Discrete and Impulsive Systems}}\\[-1.00mm]
http:monotone.uwaterloo.ca/$\sim$journal} $~$ \\ [.2in]


\begin{center}
{\large\bf ON STABILITY OF GENERALIZED CAUCHY-TYPE PROBLEM}

\vskip.20in

Sandeep P Bhairat\\[2mm]
{\footnotesize
Department of Mathematics,\\
Institute of Chemical Technology, Mumbai--400 019, (M.S) India.\\[5pt]
}
\end{center}

{\footnotesize
\noindent
{\bf Abstract.} In this paper, we study the stability of solution of initial value problem for fractional differential equation involving generalized Katugampola derivative. Pachpatte inequality is used as handy tool to obtain our result.\\[3pt]
}

\vskip.2in

\section{Introduction}

\noindent Nowadays, the subject of fractional calculus attracted great attention of many researchers and emerged as an advancement in applied mathematics. In last three decades, fractional calculus found useful for capturing naturally arising complex phenomena. The theory of arbitrary order achieved a new height in the description of properties of viscoelastic materials and memory mechanism \cite{hr2,kst}, also see \cite{bt1,fm,xu}. In recent years there has been a considerable interest in qualitative properties of fractional differential equations by using numerous operators and variety of techniques, see \cite{as},\cite{sp2}-\cite{kmf},\cite{nehh,mdk,ke,oo}.

A new fractional derivative called generalised Katugampola fractional derivative, which has unified approach, generalizes ten existing fractional derivatives (see Definition \ref{d6} below). Existence and uniqueness of solution of fractional differential equations (FDE) involving this operator are given in (Section 5, \cite{oo}). The associated fractional integral operator called Katugampola fractional integral was introduced by U Katugampola in \cite{ki} which interpolates between Riemann-Liouville and Hadamard fractional integrals.

In \cite{db5} authors have considered the initial value problem (IVP)
\begin{equation}\label{p1}
\begin{cases}
&{^{\rho}D_{a+}^{\alpha,\beta}x(t)}=f(t,x(t)),\qquad 0<\alpha<1,0\leq\beta\leq1,\rho>0,\\
&{^{\rho}I_{a+}^{1-\gamma}x(a+)}=b,\qquad b\in\mathbb{R}.
\end{cases}
\end{equation}
for FDE. The IVP \eqref{p1} is equivalent to the following integral equation
\begin{equation}\label{i1}
x(t)=\frac{b}{\Gamma(\gamma)}{\bigg(\frac{t^\rho-a^\rho}{\rho}\bigg)}^{\gamma-1}+\int_{a}^{t}s^{\rho-1}
{\bigg(\frac{t^\rho-s^\rho}{\rho}\bigg)}^{\alpha-1}\frac{f(s,x(s))}{\Gamma(\alpha)}ds.
\end{equation}
The existence and stability results are obtained using fixed point theory.

The aim of the present paper is to study the stability of generalised Cauchy-type problem involving generalized Katugampola derivative
\begin{equation}\label{1}
\begin{cases}
&{^{\rho}D_{a+}^{\alpha,\beta}x(t)}=f(t,x(t)),\qquad 0<\alpha<1,0\leq\beta\leq1,\rho>0,\\
&{\big(\frac{t^\rho-a^\rho}{\rho}\big)}^{(1-\beta)(1-\alpha)}{x(t)\big|}_{t=a}=b,\qquad b\in\mathbb{R}\backslash\{0\}.
\end{cases}
\end{equation}
Clearly, the IVP \eqref{1} is equivalent to the integral equation
\begin{equation}\label{i1}
x(t)=\frac{b}{\Gamma(\gamma)}{\bigg(\frac{t^\rho-a^\rho}{\rho}\bigg)}^{\gamma-1}+\int_{a}^{t}s^{\rho-1}
{\bigg(\frac{t^\rho-s^\rho}{\rho}\bigg)}^{\alpha-1}\frac{f(s,x(s))}{\Gamma(\alpha)}ds.
\end{equation}

The remaining paper is arranged as follows: in Section 2, we recall the preliminary facts useful for further discussion. In Section 3, we state and prove our main results. Pachpatte inequality is the main ingredient.

\section{Preliminaries}
\noindent Let us consider some definitions and basic lemmas herein.
\begin{definition}\label{d1}\cite{kst}
The space $X_{c}^{p}(a,b)\,(c\in\mathbb{R},p\geq1)$ consists of those real-valued Lebesgue measurable functions $g$ on $(a,b)$ for which ${\|g\|}_{X_{c}^{p}}<\infty,$ where
\begin{gather*}
{\|g\|}_{X_{c}^{p}}={\bigg(\int_{a}^{b}{|t^cg(t)|}^{p}\frac{dt}{t}\bigg)}^{\frac{1}{p}},\quad p\geq1,\,\,c\in\mathbb{R}\\
{\|g\|}_{X_{c}^{p=\infty}}=\text{ess sup}_{a\leq t\leq b}|t^cg(t)|,\quad c\in\mathbb{R}.
\end{gather*}
In particular, when $c=\frac{1}{p},$ we see that $X_{{1}/{p}}^{c}(a,b)=L_p(a,b).$
\end{definition}
\begin{definition}\label{d2}\cite{oo}
Let $\Omega=[a,b]$ be a finite interval on $\mathbb{R}^{+}$ and $a,\rho>0$, $0\leq\gamma<1.$
Denote by $C[a,b]$ a space of continuous functions $g$ on $\Omega$ with the norm
\begin{gather*}
{\|g\|}_{C}=\max_{t\in\Omega}|g(t)|.
\end{gather*}
The weighted space $C_{\gamma,\rho}[a,b]$ of functions $g$ on $(a,b]$ is defined by
\begin{gather}\label{space}
C_{\gamma,\rho}[a,b]=\bigg\{g:(a,b]\to\mathbb{R}:{\bigg(\frac{t^\rho-a^\rho}{\rho}\bigg)}^{\gamma}g(t)\in{C[a,b]}\bigg\},\quad0\leq\gamma<1
\end{gather}
with the norm
\begin{gather*}
{\|g\|}_{C_{\gamma,\rho}}={\bigg\|{\bigg(\frac{t^\rho-a^\rho}{\rho}\bigg)}^{\gamma}g(t)\bigg\|}_{C}=\max_{t\in\Omega}\bigg|{\bigg(\frac{t^\rho-a^\rho}{\rho}\bigg)}^{\gamma}g(t)\bigg|,
\end{gather*}
and $C_{0,\rho}[a,b]=C[a,b].$
\end{definition}
\begin{definition}\label{d3}\cite{oo}
Let $\delta_\rho=\big(t^{\rho-1}\frac{d}{dt}\big),\,\Omega=[a,b]\,(0<a<b<\infty)$ and $\rho>0,\,0\leq\gamma<1.$ For $n\in\N,$ denote $C_{\delta_\rho,\gamma}^{n}[a,b]-$ the Banach space of functions $g$ which are continuously differentiable, with $\delta_\rho,$ on $[a,b]$ upto $(n-1)$ order and have the derivative $\delta_\rho^ng$ on $(a,b]$ such that $\delta_\rho^ng\in{C_{\gamma,\rho}[a,b]},$
\begin{gather*}
C_{\delta_\rho,\gamma}^{n}[a,b]=\big\{g:[a,b]\to{\mathbb{R}}|\delta_\rho^kg\in{C[a,b]}\,\,\,\text{for}\,\,\,
0\leq{k}\leq{n-1},\\
\,\,\,\delta_\rho^ng\,\,\, \text{exists and}\,\,\,\delta_\rho^ng\in{C_{\gamma,\rho}[a,b]}\big\}
\end{gather*}
with the norm
\begin{gather*}
{\|g\|}_{C_{\delta_\rho,\gamma}^{n}}=\sum_{k=0}^{n-1}{\|\delta_\rho^kg\|}_{C}+{\|\delta_\rho^ng\|}_{C_{\gamma,\rho}},\quad
{\|g\|}_{C_{\delta_\rho}^{n}}=\sum_{k=0}^{n}\max_{t\in\Omega}|\delta_\rho^kg(t)|.
\end{gather*}
\end{definition}
\begin{definition}\label{d4}\cite{kd}
Let $g\in{X_{c}^{p}(a,b)},$ where $X_{c}^{p}$ is as in Definition \ref{d1} and $\alpha>0$. The left-sided Katugampola fractional integral $^{\rho}I_{a+}^{\alpha}$  of order $\alpha$ is defined by
\begin{gather}\label{kil}
^{\rho}I_{a+}^{\alpha}g(t)=\int_{a}^{t}s^{\rho-1}{\bigg(\frac{t^\rho-s^\rho}{\rho}\bigg)}^{\alpha-1}\frac{g(s)}{\Gamma(\alpha)}ds,\quad t>a.
\end{gather}
\end{definition}
\begin{definition}\label{d5}\cite{kd}
Let $\alpha\in{\mathbb{R}^{+}{\setminus}\mathbb{N}}$ and $n=[\alpha]+1,$ where $[\alpha]$ is integer part of $\alpha$ and $\rho>0.$ The left-sided Katugampola fractional derivative $^{\rho}D_{a+}^{\alpha}$ is defined by
\begin{align}\label{kdl}
^{\rho}D_{a+}^{\alpha}g(t)&=\delta_\rho^n (^{\rho}I_{a+}^{n-\alpha}g(s))(t).
\end{align}
\end{definition}
\begin{definition}\cite{oo}\label{d6}
The left-sided generalized Katugampola fractional derivative $^{\rho}D_{a+}^{\alpha,\beta}$ of order $0<\alpha<1$ and type $0\leq\beta\leq1$ is defined by
\begin{gather}\label{gkl}
{(^{\rho}D_{a+}^{\alpha,\beta}g)}(t)={({^{\rho}I_{a+}^{\beta(1-\alpha)}}\delta_\rho{{^{\rho}I_{a+}^{(1-\beta)(1-\alpha)}}}g)}(t),
\end{gather}
for the functions for which right-hand side expression exists and $\rho>0.$
\end{definition}
\begin{lemma}\label{l10}\cite{oo}
If $\alpha>0$ and $0<\gamma\leq1,$ then ${^\rho{I}_{a+}^{\alpha}}$ is bounded from ${C_{1-\gamma,\rho}[a,b]}$ into ${C_{1-\gamma,\rho}[a,b]}.$
\end{lemma}
\begin{lemma}\label{l4}\cite{kmf}
For nonnegative $a_i,\,i=1,\cdots,k,$
\begin{equation}\label{kl1}
{\bigg(\sum_{i=1}^{k}a_i\bigg)}^{p}\leq{k^{p-1}\sum_{i=1}^{k}{a_i}^{p}},\qquad p\geq1.
\end{equation}
\end{lemma}
\begin{lemma}\label{l5}\cite{bg}[Pachpatte Lemma]
Let $a(t)$ and $b(t)$ be continuous positive functions defined on $[t_0,\infty),$ $t_0\geq0.$ Let $w:[0,\infty)\to[0,\infty)$ be a continuous monotonic nondecreasing function such that $w(0)=0$ and $w(x)>0$ for $x>0.$ If $u$ is a positive differentiable function on $[t_0,\infty)$ that satisfies
\begin{equation*}
  u'(t)\leq{a(t)w(u(t))+b(t)},\qquad t\in{[t_0,\infty)},
\end{equation*}
then we have
\begin{equation*}
u(t)\leq{G^{-1}{\bigg[G\bigg(u(t_0)+\int_{t_0}^{t}b(s)ds\bigg)+\int_{t_0}^{t}a(s)ds\bigg]}},
\end{equation*}
for the values of $t$ for which the right-hand side is well-defined, where
\begin{equation*}
G(r)=\int_{r_0}^{r}\frac{ds}{w(s)},\qquad r>r_0>0.
\end{equation*}
\end{lemma}
In order to obtain the stability of solution for generalized Cauchy-type problem \eqref{1}, we introduce the following spaces:
\begin{equation}\label{workspace}
C_{1-\gamma,\rho}^{\alpha,\beta}[a,b]=\{g\in{C_{1-\gamma,\rho}[a,b]}:{^\rho{D}_{a+}^{\alpha,\beta}}g\in{C_{1-\gamma,\rho}[a,b]}\},
\end{equation}
and
\begin{gather*}
C_{1-\gamma,\rho}^{\gamma}[a,b]=\{g\in{C_{1-\gamma,\rho}[a,b]}:{^\rho{D}_{a+}^{\gamma}}g\in{C_{1-\gamma,\rho}[a,b]}\},\,\,0<\gamma\leq1.
\end{gather*}
Since ${^\rho{D}_{a+}^{\alpha,\beta}}g={^\rho{I}_{a+}^{\beta(1-\alpha)}}{^\rho{D}_{a+}^{\gamma}}g,$ we have $C_{1-\gamma,\rho}^{\gamma}[a,b]\subset{C_{1-\gamma,\rho}^{\alpha,\beta}[a,b]}$ follows from Lemma \ref{l10}.
\section{Stability of solution}
In this section, we present stability of global solution of the Cauchy-type problem \eqref{1}. Following lemma is of great importance in further discussion.
\begin{lemma}\label{l31}
If $\zeta,\vartheta,\varpi>0,$ then
\begin{gather}\label{2}
\hspace{-.3cm}{\bigg(\frac{t^\rho-a^\rho}{\rho}\bigg)}^{1-\vartheta}\int_{a}^{t}{\bigg(\frac{t^\rho-s^\rho}{\rho}\bigg)}^{\vartheta-1}
{\bigg(\frac{s^\rho-a^\rho}{\rho}\bigg)}^{\zeta-1}s^{\rho-1}e^{-\varpi{\big(\frac{s^\rho-a^\rho}{\rho}\big)}}ds\leq{C}\varpi^{-\zeta}
\end{gather}
for $t>a>0,$ where $C$ is a positive constant independent of $t.$
\end{lemma}
\begin{proof}
Denote left-hand side of inequality \eqref{2} by $I(t).$ By the change of variable $\xi=\frac{s^\rho-a^\rho}{t^\rho-a^\rho}$ we get
\begin{gather}\label{4}
I(t)={\bigg(\frac{t^\rho-a^\rho}{\rho}\bigg)}^{1-\vartheta}\int_{0}^{1}{(1-\xi)}^{\vartheta-1}\xi^{\zeta-1}
e^{-\varpi\xi{\big(\frac{t^\rho-a^\rho}{\rho}\big)}}d\xi.
\end{gather}
Observe that, for $\xi\geq1$ and $[\zeta]+1\geq\zeta,$ we have $\xi^{[\zeta]+1}\geq\xi^{\zeta}.$\newline
Since $\zeta+2\geq[\zeta]+2$ and the Gamma function is increasing in $[2,\infty),$ we have
\begin{gather*}
\Gamma(\zeta+2)\geq\Gamma([\zeta]+2)\quad\text{or}\quad\frac{1}{\Gamma([\zeta]+2)}\geq\frac{1}{\Gamma(\zeta+2)}.
\end{gather*}
Moreover we have $e^{\xi}\geq\frac{\xi^{[\zeta]+1}}{\Gamma([\zeta]+2)}$ and hence
\begin{gather}\label{5}
e^{\xi}\geq\frac{\xi^{[\zeta]+1}}{\Gamma([\zeta]+2)}\geq\frac{\xi^{\zeta}}{\Gamma([\zeta]+2)}\geq\frac{\xi^{\zeta}}{\Gamma(\zeta+2)}\quad\text{implies}\quad e^{-\xi}\leq\frac{\Gamma(\zeta+2)}{\xi^{\zeta}}.
\end{gather}
Therefore, for $0\leq\xi<\frac{1}{2}$ we obtain
\begin{gather}\label{7}
{(1-\xi)}^{\vartheta-1}\leq\max{(1,2^{1-\vartheta})}.
\end{gather}
For $\frac{1}{2}<\xi\leq1$ and $t$ such that $\varpi\xi{\big(\frac{t^\rho-a^\rho}{\rho}\big)}\geq1,$ we have
\begin{align}\label{8}
e^{-\varpi\xi{\big(\frac{t^\rho-a^\rho}{\rho}\big)}}\leq\frac{\Gamma(\zeta+2)}
{{\big(\varpi\xi{\big(\frac{t^\rho-a^\rho}{\rho}\big)}\big)}^{\zeta}}\leq\frac{\varpi^{-\zeta}}{\xi}{\Gamma(\zeta+2)}
\leq2\varpi^{-\zeta}\Gamma(\zeta+2).
\end{align}
Thus, using inequalities \eqref{5}-\eqref{8}, we obtain
\begin{align*}
{{\bigg(\frac{t^\rho-a^\rho}{\rho}\bigg)}^{\zeta}}&{(1-\xi)}^{\vartheta-1}\xi^{\zeta-1}e^{-{\varpi\xi{\big(\frac{t^\rho-a^\rho}{\rho}\big)}}}\nonumber\\
&\leq\begin{cases}
\max{(1,2^{1-\vartheta})}{{\big(\frac{t^\rho-a^\rho}{\rho}\big)}^{\zeta}}\xi^{\zeta-1} e^{-{\varpi\xi{\big(\frac{t^\rho-a^\rho}{\rho}\big)}}},& \,\, 0\leq\xi<\frac{1}{2} \\
2{(1-\xi)}^{\vartheta-1}\Gamma(\zeta+2)\varpi^{-\zeta}, & \,\, \frac{1}{2}<\xi\leq1.
\end{cases}
\end{align*}
As a consequence,
\begin{align}\label{9}
I(t)\leq\max{(1,2^{1-\vartheta})}&{{\bigg(\frac{t^\rho-a^\rho}{\rho}\bigg)}^{\zeta}}
\int_{0}^{\frac{1}{2}}\xi^{\zeta-1} e^{-{\varpi\xi{\big(\frac{t^\rho-a^\rho}{\rho}\big)}}}d\xi\nonumber\\
&~~~~+2\varpi^{-\zeta}\Gamma(\zeta+2)\int_{\frac{1}{2}}^{1}{(1-\xi)}^{\vartheta-1}d\xi.
\end{align}
A substitution $u={\varpi\xi{\big(\frac{t^\rho-a^\rho}{\rho}\big)}}$ yields that
\begin{align*}
I(t)\leq\max{(1,2^{1-\vartheta})}&{{\bigg(\frac{t^\rho-a^\rho}{\rho}\bigg)}^{\zeta}}
\int_{0}^{\infty}
{\bigg(\frac{u}{{\varpi{\big(\frac{t^\rho-a^\rho}{\rho}\big)}}}\bigg)}^{\zeta-1}
\frac{e^{-u}}{{\varpi{\big(\frac{t^\rho-a^\rho}{\rho}\big)}}}du\\
&+2\varpi^{-\zeta}\Gamma(\zeta+2){\bigg[\frac{{-(1-\xi)}^{\vartheta}}{\vartheta}\bigg]{\bigg|}_{\xi=\frac{1}{2}}^{1}}.
\end{align*}
This gives
\begin{equation}\label{10}
I(t)\leq\max{(1,2^{1-\vartheta})}\varpi^{-\zeta}\Gamma(\zeta)+\frac{2^{1-\vartheta}\varpi^{-\zeta}\Gamma(\zeta+2)}{\vartheta},
\end{equation}
which results in
\begin{equation*}
I(t)\leq\max{\{1,2^{1-\vartheta}\}}\varpi^{-\zeta}\Gamma(\zeta){\bigg(1+\frac{\zeta(\zeta+1)}{\vartheta}\bigg)}.
\end{equation*}
For $0<\eta<1,\,e^{\eta}\geq1$, therefore $\Gamma(\zeta+2)e^{\eta}\geq1\geq\eta^{\zeta}$
holds and for $t$ such that $0<{\varpi\xi{\big(\frac{t^\rho-a^\rho}{\rho}\big)}}<1$ one can proceed in a similar way to conclude the lemma with
$C=\max{\{1,2^{1-\vartheta}\}}\Gamma(\zeta){\big(1+\frac{\zeta(\zeta+1)}{\vartheta}\big)}.$
\end{proof}
Now we are ready to present our main stability result by using Lemma \ref{l31}. We introduce the following hypotheses.
\begin{description}
  \item[(H1)] $f(\cdot,x(\cdot))\in{C_{1-\gamma,\rho}^{\beta(1-\alpha)}(a,\infty)}$ for any $x\in{C_{1-\gamma,\rho}[a,\infty)}$ such that
\begin{equation}\label{13}
|f(t,x(t))|\leq{{\bigg(\frac{t^\rho-a^\rho}{\rho}\bigg)}^{\mu}e^{-\sigma\rho{{\big(\frac{t^\rho-a^\rho}{\rho}\big)}}}\phi(t){|x(t)|}^m},
\quad t>a>0,\mu\geq0,
\end{equation}
where $m\in{\N}\backslash \{1\}$ and $\phi$ is nonnegative continuous function on $[a,\infty).$
  \item[(H2)] For some $q>\frac{1}{\alpha}$ and
  $\phi(t){\big(\frac{t^\rho-a^\rho}{\rho}\big)}^{-m\beta(1-\alpha)}\in{L^{q}(a,\infty)}$ such that
\begin{equation*}
{\bigg({\bigg\|\phi(t)\bigg\|}_{q}\bigg)}^{m-1}{\bigg\|{\bigg(\frac{t^\rho-a^\rho}{\rho}\bigg)}^{-m\beta(1-\alpha)}\phi(t)\bigg\|}_{q}<K,
\end{equation*}
where
\begin{align*}
  K&={\bigg(\frac{{{(\Gamma(\alpha))}^{mq}}a^m}{{{|b|}^{mq(m-1)}}(m-1)2^{q(m+\alpha-1)-1}}\bigg)}^{1/q}{\bigg(\frac{{(p\sigma\rho)}^{\lambda_1m}}{{\Gamma(\lambda_1)}^m{(1+\frac{\lambda_1}{\lambda_2})}^m}\bigg)}^{1/p},\\
  \lambda_1&=1+p[\mu-(1-\gamma)m],\quad \lambda_2=1+p(\alpha-1),\,\,{\mu>(m-1)(1-\gamma)},
\end{align*}
and $p$ is the conjugate exponent of $q,$ i.e. $pq=p+q.$
\end{description}
\begin{theorem}
Let $0<\alpha<1,\,0\leq\beta\leq1,$ and $\gamma=\alpha+\beta(1-\alpha).$ Suppose that $f$ satisfies {(H1)} and $\phi$ satisfies {(H2)}. Then, for any solution of Cauchy-type problem \eqref{1}, there exists a positive constant $C$ such that
\begin{equation*}
{|x(t)|}\leq{C}{\bigg(\frac{t^\rho-a^\rho}{\rho}\bigg)}^{\gamma-1},\quad t>a>0.
\end{equation*}
\end{theorem}
\begin{proof}
Cauchy-type problem \eqref{1} is equivalent to the following Volterra integral equation \eqref{i1}.
Multiply both sides of \eqref{i1} by ${\big(\frac{{t}^\rho-a^\rho}{\rho}\big)}^{1-\gamma}$ and using inequality \eqref{13}, we get
\begin{align}\label{15}
{\bigg(\frac{{t}^\rho-a^\rho}{\rho}\bigg)}^{1-\gamma}|x(t)|&\leq |b|+\frac{{\big(\frac{{t}^\rho-a^\rho}{\rho}\big)}^{1-\gamma}}{\Gamma(\alpha)}
\int_{a}^{t}s^{\rho-1}{{\bigg(\frac{{t}^\rho-s^\rho}{\rho}\bigg)}}^{\alpha-1}
{{\bigg(\frac{s^\rho-a^\rho}{\rho}\bigg)}}^{\mu}\nonumber\\
&\hspace{3cm}\times e^{-\sigma\rho\big(\frac{s^\rho-a^\rho}{\rho}\big)}
\phi(s){|x(s)|}^{m}ds.
\end{align}
Let us denote the left-hand side of \eqref{15} by $y(t)$. Then inserting the terms ${\big(\frac{s^\rho-a^\rho}{\rho}\big)}^{m(1-\gamma)}{\big(\frac{s^\rho-a^\rho}{\rho}\big)}^{-m(1-\gamma)}$
inside the integral gives
\begin{align}\label{16}
y(t)\leq |b|+&\frac{{\big(\frac{{t}^\rho-a^\rho}{\rho}\big)}^{1-\gamma}}{\Gamma(\alpha)}
\int_{a}^{t}s^{\rho-1}{{\bigg(\frac{{t}^\rho-s^\rho}{\rho}\bigg)}}^{\alpha-1}
{{\bigg(\frac{s^\rho-a^\rho}{\rho}\bigg)}}^{\mu-(1-\gamma)m}\nonumber\\
&\hspace{2cm}\times e^{-\sigma\rho\big(\frac{s^\rho-a^\rho}{\rho}\big)}
\phi(s){y^{m}(s)}ds,\quad t>a>0.
\end{align}
Applying Holder inequality, 
we have
\begin{align*}
\int_{a}^{t}&s^{\rho-1}{{\bigg(\frac{{t}^\rho-s^\rho}{\rho}\bigg)}}^{\alpha-1}{{\bigg(\frac{s^\rho-a^\rho}
{\rho}\bigg)}}^{\mu-(1-\gamma)m}e^{-\sigma\rho\big(\frac{s^\rho-a^\rho}{\rho}\big)}\phi(s){y^{m}(s)}ds\\
\leq&{\bigg[\int_{a}^{t}s^{p(\rho-1)}{{\bigg(\frac{{t}^\rho-s^\rho}{\rho}\bigg)}}^{p(\alpha-1)}{{\bigg(\frac{s^\rho-a^\rho}
{\rho}\bigg)}}^{p(\mu-(1-\gamma)m)}e^{-p\sigma\rho\big(\frac{s^\rho-a^\rho}{\rho}\big)}ds\bigg]}^{1/p}\\
&\hspace{2cm}\times{\bigg[\int_{a}^{t}\phi^{q}(s){y^{qm}(s)}ds\bigg]}^{1/q},\qquad\qquad\qquad t>a>0.
\end{align*}
Since $x\in{C_{1-\gamma,\rho}(a,\infty)}$ and $\phi$ satisfies assumption ${(H2)},$ the second integral on the right hand side is finite for each fixed $t.$

Again by hypotheses ${(H2)},$ we have $\lambda_1>0,\lambda_2>0,(p\sigma\rho)>0.$ Thus, $\lambda_1-1=p[\mu-(1-\gamma)m]>0\,\,,\lambda_2-1=(\alpha-1)p>0.$ Thanks to Lemma \ref{l31}, 
we obtain
\begin{align}\label{18}
\int_{a}^{t}s^{\rho-1}{{\bigg(\frac{{t}^\rho-s^\rho}{\rho}\bigg)}}^{\alpha-1}&{{\bigg(\frac{s^\rho-a^\rho}
{\rho}\bigg)}}^{\mu-(1-\gamma)m}e^{-\sigma\rho\big(\frac{s^\rho-a^\rho}{\rho}\big)}\phi(s){y^{m}(s)}ds\nonumber\\
&\leq{C_1{{\bigg(\frac{{t}^\rho-a^\rho}{\rho}\bigg)}}^{\alpha-1}}{\bigg[\int_{a}^{t}\phi^{q}(s){y^{qm}(s)}ds\bigg]}^{1/q},
\end{align}
with $C_1={[2^{(\alpha-1)p}\Gamma(\lambda_1)(1+\frac{\lambda_1(\lambda_1+1)}{\lambda_2}){(p\sigma\rho)}^{-\lambda_1}]}^{1/p}.$
Linking \eqref{16} and \eqref{18} we obtain
\begin{equation}\label{19}
y(t)\leq|b|+{{\hat{C}}_1{{\bigg(\frac{{t}^\rho-a^\rho}{\rho}\bigg)}}^{-\beta(\alpha-1)}}
{\bigg(\int_{a}^{t}\phi^{q}(s){y^{qm}(s)}ds\bigg)}^{1/q},\,\, t>a>0,
\end{equation}
for ${\hat{C}}_1=\frac{C_1}{\Gamma(\alpha)}.$ Multiply to both sides of \eqref{19} by ${\big(\frac{{t}^\rho-a^\rho}{\rho}\big)}^{\beta(\alpha-1)},$ we obtain
\begin{equation}\label{20}
{\bigg(\frac{{t}^\rho-a^\rho}{\rho}\bigg)}^{\beta(\alpha-1)}y(t)
\leq|b|{\bigg(\frac{{t}^\rho-a^\rho}{\rho}\bigg)}^{\beta(\alpha-1)}+{\hat{C}}_1
{\bigg(\int_{a}^{t}\phi^{q}(s){y^{qm}(s)}ds\bigg)}^{1/q}.
\end{equation}
Denote by $z(t)$ the left-hand side of \eqref{20}. Insert the term ${\big(\frac{s^\rho-a^\rho}{\rho}\big)}^{-qm\beta(\alpha-1)}$ ${\big(\frac{s^\rho-a^\rho}{\rho}\big)}^{qm\beta(\alpha-1)}$
inside the integral on the right-hand side of \eqref{20} gives
\begin{equation}\label{21}
z(t)\leq|b|{\bigg(\frac{{t}^\rho-a^\rho}{\rho}\bigg)}^{\beta(\alpha-1)}+{\hat{C}}_1
{\bigg[\int_{a}^{t}\phi^{q}(s){\bigg(\frac{s^\rho-a^\rho}{\rho}\bigg)}^{-qm\beta(\alpha-1)}{z^{qm}(s)}ds\bigg]}^{1/q}.
\end{equation}
Raising both sides of \eqref{21} to the power $q$, we get
\begin{equation}\label{22}
z^q(t)\leq2^{q-1}{\bigg[{|b|}^{q}{\bigg(\frac{{t}^\rho-a^\rho}{\rho}\bigg)}^{q\beta(\alpha-1)}+{{\hat{C}}_1}^q
\int_{a}^{t}\phi^{q}(s){\bigg(\frac{s^\rho-a^\rho}{\rho}\bigg)}^{-qm\beta(\alpha-1)}{z^{qm}(s)}ds\bigg]}.
\end{equation}
Set
\begin{equation}\label{23}
w(t)={{\hat{C}}_1}^q\int_{a}^{t}\phi^{q}(s){\bigg(\frac{s^\rho-a^\rho}{\rho}\bigg)}^{-qm\beta(\alpha-1)}{z^{qm}(s)}ds,\quad t>a>0.
\end{equation}
Then, by the continuity of $z(t)$ and assumption ${(H2)},$ the integrand is summable. Clearly $w(a)=0,$ and by differentiation
\begin{equation}\label{24}
w'(t)={{\hat{C}}_1}^q\phi^{q}(t){\bigg(\frac{t^\rho-a^\rho}{\rho}\bigg)}^{-qm\beta(\alpha-1)}{z^{qm}(t)},\quad t>a>0.
\end{equation}
Moreover, $\phi,z$ and right-hand side of \eqref{24} are nonnegative, $w$ is a nonnegative continuous and nondecreasing function in $[a,\infty).$

Further, we estimate the right-hand side of \eqref{24} in terms of $w(t).$ From \eqref{22} and \eqref{23}, we obtain
\begin{equation*}
z^q(t)\leq2^{q-1}{\bigg[{|b|}^{q}{\bigg(\frac{{t}^\rho-a^\rho}{\rho}\bigg)}^{q\beta(\alpha-1)}+w(t)\bigg]}.
\end{equation*}
Raising both sides to the power $m$ and using Lemma \ref{l4}, we get
\begin{equation}\label{25}
z^{qm}(t)\leq2^{mq-1}{\bigg[{|b|}^{mq}{\bigg(\frac{{t}^\rho-a^\rho}{\rho}\bigg)}^{mq\beta(\alpha-1)}+w^m(t)\bigg]}.
\end{equation}
Substituting \eqref{25} into \eqref{24} yields
\begin{align}\label{26}
w'(t)&\leq2^{mq-1}{\hat{C}_1}^{q}{\phi}^q(t){\bigg(\frac{{t}^\rho-a^\rho}{\rho}\bigg)}^{-mq\beta(\alpha-1)}\nonumber\\
&\hspace{2cm}\times{\bigg[{|b|}^{mq}{\bigg(\frac{{t}^\rho-a^\rho}{\rho}\bigg)}^{mq\beta(\alpha-1)}+w^m(t)\bigg]}\nonumber\\
&\leq2^{mq-1}{|b|}^{mq}{\hat{C}_1}^{q}{\phi}^q(t)\nonumber\\
&\hspace{2cm}+2^{mq-1}{\hat{C}_1}^{q}
{\bigg(\frac{{t}^\rho-a^\rho}{\rho}\bigg)}^{-mq\beta(\alpha-1)}{\phi}^q(t)w^m(t).
\end{align}
Applying Pachpatte Lemma (Lemma \ref{l5} with $w(x)=x^m$) we infer that
\begin{align}\label{27}
w(t)\leq{G^{-1}}&\bigg[G\bigg(w(a)+2^{mq-1}{|b|}^{mq}{\hat{C}_1}^q\int_{a}^{t}\phi^q(s)ds\bigg)\nonumber\\
&+2^{mq-1}{\hat{C}_1}^q\int_{a}^{t}{\bigg(\frac{s^\rho-a^\rho}{\rho}\bigg)}^{-mq\beta(\alpha-1)}\phi^q(s)ds\bigg].
\end{align}
Set
\begin{gather*}
l(t)=2^{mq-1}{|b|}^{mq}{\hat{C}_1}^q\int_{a}^{t}\phi^q(s)ds,\\
k(t)=2^{mq-1}{\hat{C}_1}^q\int_{a}^{t}{\bigg(\frac{s^\rho-a^\rho}{\rho}\bigg)}^{-mq\beta(\alpha-1)}\phi^q(s)ds,
\end{gather*}
then inequality \eqref{27} becomes
\begin{equation}\label{29}
w(t)\leq{G^{-1}}\big[G\big(l(t)\big)+k(t)\big],
\end{equation}
where we have used the fact $w(a)=0.$ Here $G(r)=\int_{r_0}^{r}\frac{ds}{s^m}$, $r>0,\,r_0>0,$
\begin{equation*}
i.e.\,\,\,G(r)=\frac{r^{1-m}}{1-m}-\frac{{r_0}^{1-m}}{1-m}\quad\text{and}\quad G^{-1}(y)=[{r_0}^{1-m}-(m-1)y]^{-1/(m-1)}.
\end{equation*}
The inequality \eqref{29} reads
\begin{align}\label{30}
w(t)&\leq G^{-1}{\bigg[\frac{{l(t)}^{1-m}}{1-m}-\frac{{l(t_0)}^{1-m}}{1-m}+k(t)\bigg]}\nonumber\\
&\leq{\bigg[{l(t_0)}^{1-m}-(m-1)\bigg(\frac{{l(t)}^{1-m}}{1-m}-\frac{{l(t_0)}^{1-m}}{1-m}+k(t)\bigg)\bigg]}^{-\frac{1}{m-1}}\nonumber\\
&\leq{[{l(t)}^{1-m}-(m-1)k(t)]}^{-\frac{1}{m-1}},
\end{align}
as long as ${l(t)}^{m-1}k(t)<{\frac{1}{m-1}}.$ In particular, if
\begin{equation*}
{\bigg(\int_{a}^{t}\phi^q(s)ds\bigg)}^{m-1}{\bigg[\int_{a}^{t}\phi^q(s){\bigg(\frac{s^\rho-a^\rho}{\rho}\bigg)}^{-mq\beta(1-\alpha)}ds\bigg]}<K/2
\end{equation*}
then $w(t)\leq{K_1}$ for some positive constant $K_1$ for all $t>a>0$ and thus from \eqref{21}, we find that
\begin{equation*}
  z(t)\leq|b|{\bigg(\frac{t^\rho-a^\rho}{\rho}\bigg)}^{\beta(1-\alpha)}+{K_1}^{1/q},
\end{equation*}
and then
\begin{equation*}
y(t)\leq|b|+{K_1}^{1/q}{\bigg(\frac{t^\rho-a^\rho}{\rho}\bigg)}^{-\beta(1-\alpha)}\leq{C},\quad t\geq{t_0}>a>0,
\end{equation*}
for some positive constant $C.$ This yields that
\begin{equation*}
  |x(t)|\leq{C}{{\bigg(\frac{t^\rho-a^\rho}{\rho}\bigg)}^{\gamma-1}}\qquad\text{for}\quad t\geq{t_0}>a>0.
\end{equation*}
\end{proof}



\end{document}